\documentclass{ifacconf}
\usepackage{graphicx}     
\usepackage{natbib}   
\usepackage{graphicx} 
\usepackage{amsmath}
\usepackage{eurosym}
\usepackage{textcomp}
\usepackage{ulem}
\usepackage{subcaption}
\usepackage[utf8]{inputenc}
\usepackage{pgfplots}
\usepackage{xurl}
\usepackage{bm}
\DeclareUnicodeCharacter{2212}{−}
\usepgfplotslibrary{groupplots,dateplot}
\usetikzlibrary{patterns,shapes.arrows}
\pgfplotsset{compat=newest}

\begin{document}
\begin{frontmatter}
\title{Backcasting Policies in Transport Systems as an Optimal Control Problem : An Example with Electric Vehicle Purchase Incentives}

\author{Vinith Lakshmanan} 
\author{Xavier Guichet} 
\author{Antonio Sciarretta}

\address{IFP Energies Nouvelles, France (e-mail: vinith-kumar.lakshmanan@ifpen.fr, xavier.guichet@ifpen.fr, antonio.sciarretta@ifpen.fr).}

\begin{abstract}                
This study represents a first attempt to build a backcasting methodology to identify the optimal policy roadmaps in transport systems. Specifically, it considers a passenger car fleet subsystem, modelling its evolution and greenhouse gas emissions. The policy decision under consideration is the monetary incentive to the purchase of electric vehicles. This process is cast as an optimal control problem with the objective to minimize the total budget of the state and reach a desired CO$_2$ target. A case study applied to Metropolitan France is presented to illustrate the approach. Additionally, alternative policy scenarios are also analyzed.   
\end{abstract}


\end{frontmatter}

\section{Introduction}

The European Union's (EU) goal of carbon neutrality by 2050 requires reducing transportation sector emissions by 90\% compared to 1990 levels. Transport alone is set to make up nearly half of Europe's greenhouse gas (GHG) emissions in 2030. Thus the European Commission has adopted a set of proposals to make the EU's policies fit for reducing net greenhouse gas emissions by at least 55\% by 2030 (\cite{GreenDeal}).
Similar targets are being implemented in U.S. and other regions (\cite{EPA}). 

Governance, policies, and incentives ("decisions") play an important role in shaping transport systems of the future, and influence the development and implementation of the various technologies and modes of transport. It is therefore important to study how decisions could be best used to govern transport systems in the desired direction of decarbonization. 


To find the best policy roadmaps for desired targets, the traditional approach consists in designing prospective scenarios, and testing them using simulation. Once the impacts are forecasted for each scenario, conclusions can be drawn on which decisions are the most effective. With this approach, the choice is limited to the scenarios designed, which may represent a tiny subset of all possibilities if multiple concurrent decisions are considered. Moreover, even when a small set of decisions are to be taken, the optimum might not be achieved since only the designed scenarios are evaluated.

To overcome these limitations, a backcasting paradigm is supported in this work. In this approach, desired targets are set by the decision makers at a certain time horizon, then the optimal combinations of policies to achieve these targets are calculated as a function of time (``backcasted"). 
In this way, the aprioristic choice of scenarios is replaced by a full dynamic optimization process that can explore among all combinations possible.

The backcasting paradigm has been introduced since the last century (\cite{Robinson, Bibri}). It has been mainly deployed in qualitative terms (\cite{Levitate}) or quantitatively with some static optimization procedure (\cite{Gomi, Ashina}).
However, this process can be more effectively cast as an optimal control problem, with a suitable definition of an objective function, an horizon, local and terminal constraints, etc.


Clearly, since future impacts have to be predicted, the new backcasting paradigm is still based on a simulation model. This model must be able to describe transport as a system, with manipulable inputs, exogenous inputs or disturbances, outputs, and states. The manipulable inputs are represented by the decisions to optimize, which may concern local authorities, state government, EU, or even private companies. 
The exogenous inputs represent the influence of other, related systems such as the energy, urban, economic, demographic ones, which cannot be modeled within the transport system alone.
The outputs represent the impacts targeted or the constraints to impose to the backcasting process. Finally, the states are the dynamics associated with the internal variables. 

In this paper, we illustrate the backcasting paradigm, applied to the transport sector, by considering a specific subsystem with a single decision variable. The subsystem considered describes the evolution of the passenger car fleet within a certain region and its impacts on the GHG emissions. The decision optimized is a monetary incentive to the purchase of electric vehicles.

The prediction of vehicle fleet composition is the subject of a large body of literature (\cite{TREMOVE, High-Tool, ITF, ADEME}). Typically, dynamic fleet models are based on the evaluation of stocks and sales of various types of vehicles per time period. Stocks change in time due to disposal of old vehicles (due to scrappage, exports, change of use, etc.), and sales of new vehicles. The latter, in turn, are induced by transport demand (vkm) and mileage, and split among the vehicle types using discrete choice models (\cite{train, benakiva}). 

The GHG emissions of a given vehicle fleet are typically assessed using emission factors. CO$_2$ emissions of light-duty vehicles are regulated in the EU. Similar regulations are about to be applied to heavy-duty vehicles as well.

Recent studies that include electric vehicles have applied a fleet model to predict the future transport emissions in France (\cite{ADEME, ITF}), Norway (\cite{thorne}), Japan (\cite{sato}), and the U.S. (\cite{woody}).

The paper is organized as follows. Section~\ref{sec:model-refinement} introduces the passenger car fleet model, followed by the optimal control problem formulation. A case study, based on the French national passenger car fleet, is presented in Sect.~\ref{sec:casestudy}. The last section draws some conclusions and proposes several research directions to extend to a more realistic system model.

\section{Passenger Car Fleet Model} \label{sec:model-refinement}
This section presents the equations of the passenger car fleet model and the formulation of the proposed backcasting approach.  The latter consists in optimising some decisions concerning the transport system  in order to achieve some defined target in greenhouse-gas (CO$_2$) emissions at year $T$. We thus treat the transport system as a system having $\textbf{u}=
u_v(t)$ as the manipulable input to be backcasted in this study, where the monetary incentive to the purchase of electric vehicles (EV) given by the state ($u\equiv u_2$, while $u_1\equiv 0$), and the CO$_2$ emissions as the targeted output. 
This system shall be represented by an aggregated dynamic model, which evaluates the transport emissions of the studied area as a function of time.

\subsection{Model}

We consider only a single area of interest and private car as the transport mode. In addition, we consider the latter's stock composed of two types of vehicles (thermal, $v=1$ and electric, $v=2$) differentiated by $A+1$ classes of ages ($a=0\ldots A$). We consider vehicle-km (vkm), a measure of transport demand $G(t)$, as an input provided by upstream models. Additionally, we also consider mileage $M(t)$ as an exogenous input.

We use one year as the time step and label the year index as $t$, starting from present until target year $T$. The passenger car fleet model can be written as follows. 

The demand of new vehicles $N$ at year $t$ is given by the ratio of the vkm demand for new vehicles, $F(t)$ and the mileage $M(t)$,
\begin{equation} \label{eqn:N}
  N(t) = \frac{F(t)}{M(t)}.
\end{equation}
The vkm demand for new vehicles is evaluated as the difference between total transport demand, $G(t)$ and those covered by the sum of the old vehicles disaggregated by vehicle type and age $O_{va}(t)$,
\begin{equation}\label{eqn:F2}
    F(t) = G(t) - M(t) \sum_{va} O_{va}(t)
\end{equation}
The latter is obtained using the age-dependant survival rate $\eta_a$ and the stock by vehicle type and age $S_{va}(t)$ at year $t$ as 
\begin{equation}\label{eqn:O2}
O_{va}(t) = \left\{ \begin{array}{ll} 
\eta_a S_{v,a-1}(t-1)  \quad \forall a=1,\ldots,A-1 \\  
\eta_A S_{v,A-1}(t-1) + \eta_A S_{vA}(t-1) 
\end{array} \right.
\end{equation}
%
The total sales at year $t$ are split among vehicle types according to 
\begin{equation} \label{eqn:Nv}
  N_v(t) = P_v(t)N(t),
\end{equation}
where $P_v$ is the share of sales by veh-type at year $t$. 

The latter is obtained from a logit expression
\begin{equation} \label{eqn:P}
  P_v(t) = \frac{e^{\mu U_v(t)}}{\sum_v e^{\mu U_v(t)} },
\end{equation}
where $U_v$ is the utility function by veh-type.

To evaluate the $P_v(t)$, we consider different technical-economic characteristics of the vehicles in its utility function, $U_v(t)$. Among the latter, we consider two classes of costs for the user: purchase costs, and operating costs (i.e., fuel/energy, maintenance, insurance costs). 
These are the main determinants for the choice of new vehicles. Another determinant is the development rate of the refilling (fuel/electricity) infrastructure which reflects its availability.
Conversely, we do not consider explicitly socio-economic determinants that depend on the single agent (like age, gender, income level, etc.) and thus are difficult to be accounted for in an aggregated model. 
Instead, we introduce an adoption coefficient to better model the penetration of new technologies (\cite{macmanus, sterman}) such as the EV. The latter is a probability density function based on the Bass model (\cite{Bass}). The expression for $U_v(t)$ is given as
\begin{equation} \label{eqn:U}
  \begin{aligned}
  U_v(t) = \left(1-c_v^A(t) \right)  \Biggl(p^P\frac{C_v^P(t)-u(t)}{\overline{C}^P(t)} + p^O \frac{C_v^O(t)}{\overline{C}^O(t)} +  \\   + p^I {(1-c_v^I(t))} \Biggr),  
  \end{aligned}
\end{equation}  
where $p$'s are tuning coefficients, $C_v^P$ is the purchase price, $C^O$ is the sum of operational costs, and $c^I$ is the rate of development of the refilling infrastructure (normalized to unity, by definition $c_1^I\equiv 1$). The average costs between the two vehicles types are given by $\overline{C}_v^P(t)$ and $\overline{C}_v^O(t)$. The prefactor that multiplies the cost-based utility is the adoption coefficient $c^A_2$ ($c^A_1\equiv 0$). Since the cost-based utility is negative (coefficients $p$'s are so), a prefactor lower than unity increases the utility of EVs proportionally to their rate of exposure.

The evolution of vehicle stock by vehicle type and age, $s_{va}$ at year $t$, is given by 
\begin{equation}
    S_{va}(t) = \left\{ \begin{array}{lll} N_v(t) & a=0 \\ O_{va}(t) & a \geq 1 \end{array} \right. .
\end{equation}

Tailpipe CO$_2$ emissions are described by simple emission factors (g/km) in this work. The latter are certainly differentiated by vehicle type and generally with vehicle age, since vehicles produced in a certain year have to comply with the emission regulations in force that year. The emissions of the stock are evaluated using the age-specific factor $\epsilon_{va}(t)$ and its annual mileage $M(t)$ as
\begin{equation}
    E(t) = \sum_{va} \epsilon_{va}(t) M(t) S_{va}(t).
\end{equation}

\subsection{OCP Formulation}

If $T$ is the time horizon, we state an optimal control problem where the cost function is the total budget for the state $I(T)$, that is, the sum of yearly products of the incentive and the number of EV sales
\begin{equation} \label{eqn:Jr}
        \min_{u(t)} I(T) = \sum_{t=1}^T  u(t) (1- P_1(t,u(t))) N(t, \textbf{x}(t-1))
\end{equation}
where the explicit form of $N$ is
\begin{equation}
\begin{aligned}
     & N(t,\textbf{x}(t-1)) = \frac{G(t)}{M(t)} - \Biggl(\sum_{v,a=1}^{A-1} \eta_a S_{v,a-1}(t-1) + \\&  + \eta_A S_{vA-1}(t-1) +  \eta_A S_{vA}(t-1) \Biggr)\;,
\end{aligned}
\end{equation}
and the state includes all partial stocks, $\textbf{x}=[S_{v0},\ldots,S_{vA}]$, $\forall v$.
Minimization of (\ref{eqn:Jr}) is subject to the terminal condition
\begin{equation} \label{eqn:ebar}
  E(T) \leq \overline{E} \;,  
\end{equation}
where $\overline{E}$ is the desired target on emissions at horizon $T$, to state equations
\begin{equation}
\begin{aligned}
S_{v0}(t) &= P_v(\textbf{u}(t)) N(t,\textbf{x}(t-1))\;, \\
S_{va}(t) &= \eta_a S_{v,a-1}(t-1), \quad \forall a = 1,\ldots,A-1 \\
S_{vA}(t) &= \eta_A S_{v,A-1}(t-1)+\eta_A S_{vA}(t-1) \;,
\end{aligned}
\end{equation}
as well as to opportune initial and boundary conditions for the control and state variables. 

For a discrete system, the Hamiltonian is generally formed as 
\begin{equation*} 
H(t) = L(t,\textbf{u}(t)) +\bm{\lambda}(t) f\left(t,\textbf{u}(t),\textbf{x}(t-1)\right), \;\; \forall t=1, \ldots, T .  
\end{equation*}
where $L$ and $\bm{\lambda}$ are the running cost and the adjoint state vector, respectively.
If are there no constraints on the control, the necessary conditions are
\begin{equation} \label{eqn:lambda}
\bm{\lambda}(t-1) =\frac{\partial H(t)}{\partial \textbf{x}(t-1)} 
\end{equation}
\begin{equation}
\frac{\partial H(t)}{\partial \textbf{u}(t)}  =0 \text { at } \textbf{u}^* .
\end{equation}

The Hamiltonian for this study case can be formed as
\begin{equation}
    \begin{aligned}
        H(t) = & u(t) (1 - P_1(u)) N(t, S_{va}(t-1)) \\
        & + \sum_{v} \lambda_{v0}(t) P_v(u) N(t,S_{va}(t-1)) \\
        & + \sum_{v,a=1}^{A-1} \lambda_{va}(t) \eta_a S_{v,a-1}(t-1) \\
        & + \sum_v \lambda_{vA}(t) \eta_A \bigl( S_{v,A-1}(t-1) + S_{vA}(t-1) \bigr)
    \end{aligned}
\end{equation}

The first-order optimality conditions yield
\begin{equation}
    \begin{aligned}
        \frac{\partial{H(t)}}{\partial{u(t)}} = & \, N(t,S_{va}(t-1))\bigl(1 - P_1 + \\& + \frac{\partial{P_1}}{\partial{u(t)}}\left(\lambda_{1,0}(t) -  \lambda_{2,0}(t) - u(t)\right) \bigr) = 0
    \end{aligned}
\end{equation}
\begin{equation}
    \begin{aligned}
        \lambda_{va}(t-1) =& \frac{\partial{H}(t)}{\partial{S_{va}(t-1)}} = -\eta_{a+1} \left(u(t)\left(1-P_1\right) + \right.\\& \left. + \lambda_{1,0}P_1 + \lambda_{2,0}(1-P_1) \right) + \lambda_{v,a+1}\eta_{a+1}, \\& \quad a=0,\ldots,A-1
    \end{aligned}
\end{equation}
\begin{equation}
\begin{aligned}
      \lambda_{vA}(t-1) =& \frac{\partial{H}(t)}{\partial{S_{vA}(t-1)}} = -\eta_A\left( u(t)(1-P_1) + \right.\\& +\lambda_{1,0}P_1 + \left. \lambda_{2,0}(1-P_1)\right) + \lambda_{vA} \eta_{A}
\end{aligned}
\end{equation}
where we have omitted the dependency on time and control for the sake of shortness.



The nonlinear system of differential equations cannot be solved analytically to obtain the optimal incentive law. Therefore, numerical procedures must be employed. The solution is obtained using the \textit{'trust-constr'} method within the Python scipy optimize package (\cite{trustCons}).

\section{Case Study} \label{sec:casestudy}
This section describes a case study applied to Metropolitan France to illustrate the backcasting approach. 
It addresses the question, \textit{what financial incentives for electric vehicles make it possible to achieve a desired level of CO$_2$ in year $T$ while minimizing public spending during those years?} 
To this end, the analysis considers a time horizon from $t_0$ = 2022 to $T$ = 2050, with a one year time step. The CO$_2$ target is set by forecasting a reference scenario in which a constant incentive (IC) of 5 k\euro~is provided for every EV purchased. As a result, $\overline{E}$ in (\ref{eqn:ebar}) is set to $E(T)$ from this scenario. The 2022 passenger car fleet, from data source~1 in Table~\ref{tab:data}, is taken as the initial value for the states $\textbf{x}$. The model inputs and parameters are described in Sect.~\ref{sec:calib}, followed by the results of the case study. Furthermore, alternative policy scenarios are analyzed in Sect.~\ref{sec:alter}. 

\subsection{Parameter Calibration and Exogenous Inputs}\label{sec:calib}

The parameters of the model are tuned using historical data using sources listed in Table~\ref{tab:data}.

\begin{table}[t!]
    \centering
    \caption{Data sources}
    \begin{tabular}{|c|c|p{5.5cm}|}
        \hline
        Index & Parameter & Web link \\
        \hline
        1 & $s_{voa}$ & \url{www.statistiques.developpement-durable.gouv.fr/parc-et-circulation-des-vehicules-routiers} \\
        2 & $\epsilon_{1,0}$ & \url{carlabelling.ademe.fr/chiffrescles/r/evolutionTauxCo2} \\
        3 & $e_v$ & \url{www.citepa.org/fr/secten} \\
        4 & $\dot{\chi}$ & \url{www.statistiques.developpement-durable.gouv.fr/immatriculation-des-vehicules-routiers}\\
        \hline
    \end{tabular}
    \label{tab:data}
\end{table}

As for the survival rate, the identification was carried out using data source~1. The latter contains historical stock of passenger vehicles $s_{voa}(\tau)$ by technology, ownership type ($o = \{$private, professional$\}$), and age  until 2022.  As a first step, we neglect the dependency of survival rate on vehicle technology, and the movement of second-hand vehicles between ownership types.  Under these assumptions, the survival rate is defined as 
\begin{equation}
    \eta_a = \frac{\sum_{vo} s_{voa}(2022)}{\sum_{vo} s_{vo,a-1}(2021)} 
\end{equation}
and is approximated as an affine function, given by
\begin{equation}
    \eta_{a} = 1.05 - 0.01 \cdot a. 
\end{equation}
The approximated $\eta_a$ corresponds to an average vehicle life of around 11 years. The survival rate exceeds one for newer vehicles, likely due to vehicles imported from neighbouring countries that are subsequently sold in France; a common practice with professional ownership type.  This is overcome by saturating the maximum of $\eta_a$ to 1.

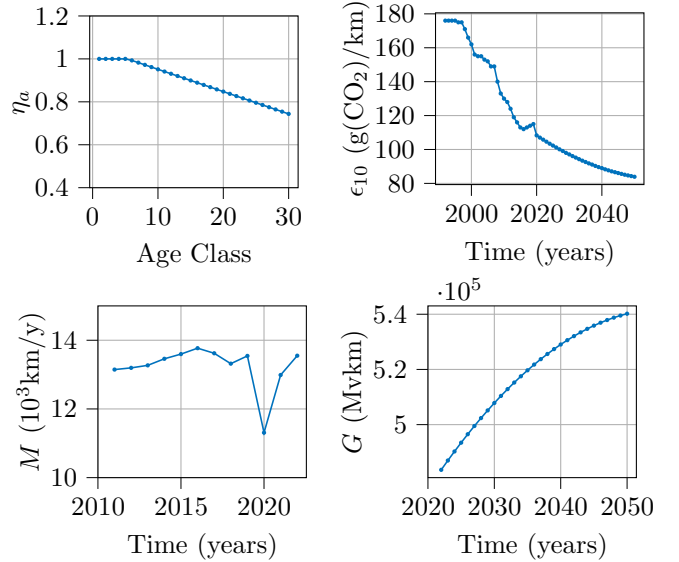
\begin{figure}[h!]
    \centering
    \begin{subfigure}[t]{0.48\columnwidth}
        \centering
        \pgfplotsset{ylabel shift = -.5em} 
\definecolor{mycolor1}{rgb}{0.00000,0.44700,0.74100}%
\definecolor{mycolor2}{rgb}{0.85000,0.32500,0.09800}%
\definecolor{mycolor3}{rgb}{0.92900,0.69400,0.12500}%
\definecolor{mycolor4}{rgb}{0.49400,0.18400,0.55600}%
\begin{tikzpicture}

\definecolor{darkgray176}{RGB}{176,176,176}
\definecolor{steelblue31119180}{RGB}{31,119,180}

\begin{axis}[
/pgf/number format/1000 sep={},
scale = 0.4,
tick align=outside,
tick pos=left,
x grid style={darkgray176},
xlabel={Age Class},
xmajorgrids,
xmin=-0.45, xmax=31.45,
xtick style={color=black},
y grid style={darkgray176},
ylabel={$\eta_a$},
ymajorgrids,
ymin=0.4, ymax=1.2,
ytick style={color=black}
]
\addplot [semithick, mycolor1, mark=*, mark size=.5, mark options={solid}]
table {%
1 1
2 1
3 1
4 1
5 1
6 0.99298212
7 0.98260041
8 0.9722187
9 0.96183699
10 0.95145528
11 0.94107357
12 0.93069186
13 0.92031015
14 0.90992844
15 0.89954673
16 0.88916502
17 0.87878331
18 0.8684016
19 0.85801989
20 0.84763818
21 0.83725647
22 0.82687476
23 0.81649305
24 0.80611134
25 0.79572963
26 0.78534792
27 0.77496621
28 0.7645845
29 0.75420279
30 0.74382108
};
\end{axis}

\end{tikzpicture}
    \end{subfigure}
    \begin{subfigure}[t]{0.48\columnwidth}
        \centering
        \pgfplotsset{ylabel shift = -.5em} 
\definecolor{mycolor1}{rgb}{0.00000,0.44700,0.74100}%
\definecolor{mycolor2}{rgb}{0.85000,0.32500,0.09800}%
\definecolor{mycolor3}{rgb}{0.92900,0.69400,0.12500}%
\definecolor{mycolor4}{rgb}{0.49400,0.18400,0.55600}%
\begin{tikzpicture}

\definecolor{darkgray176}{RGB}{176,176,176}
\definecolor{steelblue31119180}{RGB}{31,119,180}

\begin{axis}[
/pgf/number format/1000 sep={},
scale = 0.4,
tick align=outside,
tick pos=left,
x grid style={darkgray176},
xlabel={Time (years)},
xmajorgrids,
xmin=1989, xmax=2052.9,
xtick style={color=black},
y grid style={darkgray176},
ylabel={$\epsilon_{10}$ (g(CO$_2$)/km)},
ymajorgrids,
ymin=79.295, ymax=180.605,
ytick style={color=black}
]
\addplot [semithick, mycolor1, mark=*, mark size=.5, mark options={solid}]
table {%
1992 176
1993 176
1994 176
1995 176
1996 175
1997 175
1998 171
1999 166
2000 162
2001 156
2002 155
2003 155
2004 153
2005 152
2006 149
2007 149
2008 140
2009 133
2010 130
2011 128
2012 124
2013 119
2014 116
2015 113
2016 112
2017 113
2018 114
2019 115
2020 108.3
2021 107
2022 105.8
2023 104.6
2024 103.4
2025 102.3
2026 101.2
2027 100.1
2028 99.1
2029 98
2030 97.1
2031 96.1
2032 95.2
2033 94.3
2034 93.4
2035 92.6
2036 91.8
2037 91.1
2038 90.3
2039 89.6
2040 89
2041 88.3
2042 87.7
2043 87.1
2044 86.6
2045 86.1
2046 85.6
2047 85.1
2048 84.7
2049 84.3
2050 83.9
};
\end{axis}

\end{tikzpicture}
    \end{subfigure}

    \begin{subfigure}[t]{0.48\columnwidth}
        \centering
        \pgfplotsset{ylabel shift = -.5em} 
\definecolor{mycolor1}{rgb}{0.00000,0.44700,0.74100}%
\definecolor{mycolor2}{rgb}{0.85000,0.32500,0.09800}%
\definecolor{mycolor3}{rgb}{0.92900,0.69400,0.12500}%
\definecolor{mycolor4}{rgb}{0.49400,0.18400,0.55600}%

\begin{tikzpicture}

\definecolor{darkgray176}{RGB}{176,176,176}
\definecolor{steelblue31119180}{RGB}{31,119,180}

\begin{axis}[
/pgf/number format/1000 sep={},
scale=0.4,
tick align=outside,
tick pos=left,
x grid style={darkgray176},
xlabel={Time (years)},
xmajorgrids,
xmin=2010, xmax=2022.55,
xtick style={color=black},
y grid style={darkgray176},
ylabel={$M$ (10$^3$km/y)},
ymajorgrids,
ymin=10, ymax=15,
ytick style={color=black}
]
\addplot [semithick, mycolor1,  mark=*, mark size=.5, mark options={solid}]
table {%
2011 13.146
2012 13.196
2013 13.269
2014 13.460
2015 13.597
2016 13.769
2017 13.621
2018 13.318
2019 13.544
2020 11.309
2021 12.988
2022 13.550
};
\end{axis}

\end{tikzpicture}
    \end{subfigure}
    \begin{subfigure}[t]{0.48\columnwidth}
        \centering
        \pgfplotsset{ylabel shift = -.5em} 
\definecolor{mycolor1}{rgb}{0.00000,0.44700,0.74100}%
\definecolor{mycolor2}{rgb}{0.85000,0.32500,0.09800}%
\definecolor{mycolor3}{rgb}{0.92900,0.69400,0.12500}%
\definecolor{mycolor4}{rgb}{0.49400,0.18400,0.55600}%

\begin{tikzpicture}

\definecolor{darkgray176}{RGB}{176,176,176}
\definecolor{steelblue31119180}{RGB}{31,119,180}

\begin{axis}[
/pgf/number format/1000 sep={},
scale = 0.4,
tick align=outside,
tick pos=left,
x grid style={darkgray176},
xlabel={Time (years)},
xmajorgrids,
xmin=2020, xmax=2051.4,
xtick style={color=black},
y grid style={darkgray176},
ylabel={$G$ (Mvkm)},
ymajorgrids,
ymin=480791.457693442, ymax=543015.799327685,
ytick style={color=black}
]
\addplot [semithick, mycolor1, mark=*, mark size=.5, mark options={solid}]
table {%
2022 483619.836858635
2023 486997.48528258
2024 490274.587218235
2025 493451.142665625
2026 496527.151624724
2027 499502.614095582
2028 502377.530078175
2029 505151.899572478
2030 507825.722578466
2031 510398.999096238
2032 512871.72912572
2033 515243.912666961
2034 517515.549719911
2035 519686.640284596
2036 521757.184361016
2037 523727.18194917
2038 525596.633049058
2039 527365.537660657
2040 529033.895784014
2041 530601.707419081
2042 532068.972565809
2043 533435.691224346
2044 534701.863394616
2045 535867.489076597
2046 536932.568270312
2047 537897.100975762
2048 538761.087192946
2049 539524.526921864
2050 540187.420162492
};
\end{axis}

\end{tikzpicture}
    \end{subfigure}
    \caption{Model Inputs and Parameters} \label{fig:param}
\end{figure}

As for the emission factor $\epsilon_{1a}(t)$, considering only apparent tail-pipe emissions, the identification was carried using two sources. Data source 2 provides the historical trend ($\tau$ = 1995 to 2020) of average CO$_2$ emissions for newly sold ($a=0$) petrol and diesel cars. The emission factor for thermal vehicles ($v=1$) for this period is calculated as a weighted average based on the number of newly sold petrol and diesel vehicles and their respective emission factors. For vehicles sold prior to 1995, the emission factor is assumed to be at 1995 level. For the future trend ($\tau$ = 2020 to 2050), (\cite{ITF}) presents the efficiency trajectory of newly sold thermal vehicles in kWh(eq.)/100 km, projecting a 50\% improvement from 2015 to 2050. This evolution, with the initial value adjusted to align with data source 2, is converted to gCO$_2$/km and approximated using a quadratic function as
\begin{equation}
\begin{aligned}
        \epsilon_{1,0}(\tau) = &  0.01 \cdot (\tau-2020)^2 - 1.27 \cdot (\tau-2020) +\\ &+ 108.2, \quad \tau  \in[2020,2050]\;.
\end{aligned}
\end{equation}
Given a stock of thermal vehicles by age, their corresponding emission factor $\epsilon_{1a}(t)$ is obtained using the following transformation
\begin{equation}
    \epsilon_{1a} (t) = \epsilon_{1,0}(t_0+t-a).
\end{equation}
The emission factor for EVs is set to zero (i.e., $\epsilon_{2a} = 0$).
The survival rate by veh-age $\eta_a$ and emission factor of newly sold thermal vehicles $\epsilon_{1,0}(\tau)$ are shown in Fig.~\ref{fig:param}. 

For the annual mileage $M(t)$, data source 3 provides the annual CO$_2$ emissions by different vehicle types $e_v(\tau)$ in France, with $e_2(\tau) =0$. Assuming a constant annual mileage across vehicle types, $M(t)$ is calculated as 
\begin{equation}\label{eqn:M}
    M(\tau) = \frac{e_1(\tau)}{\sum_a  \sum_o s_{1oa}(\tau)\epsilon_{1,0}(\tau - a)},\; \tau \in [2011, 2022]
\end{equation}
where $s_{1a}$ represents the thermal passenger vehicle stock from data source~1, and $\epsilon_{1,0}$ is the emission factor described above. Using (\ref{eqn:M}), Fig.~\ref{fig:param} shows the annual mileage $M(\tau)$, with a dip during the COVID-19 pandemic in 2020. Overall, $M(t)$ is approximated to a constant value of 13,500~kms. Since the focus of this work is emissions reduction, we match $M(t)$ using the historical emissions data. This average is slightly higher than the one estimated in (\cite{ITF}). The latter reports a mileage of 12500~kms/year computed as a function of different use profiles specific a household's location type (rural, urban, etc.).

The parameters and determinants related to the Logit model in (\ref{eqn:U}) are considered as exogenous inputs. The assumptions and values for these inputs are directly adopted from the work conducted for the French Agency of Ecological Transistion (ADEME) using the DRIVE$^{RS}$ fleet model, see (\cite{ADEME}). The transport demand $G(t)$ (vkm), shown in Fig.~\ref{fig:param}, is also adopted from the latter.  

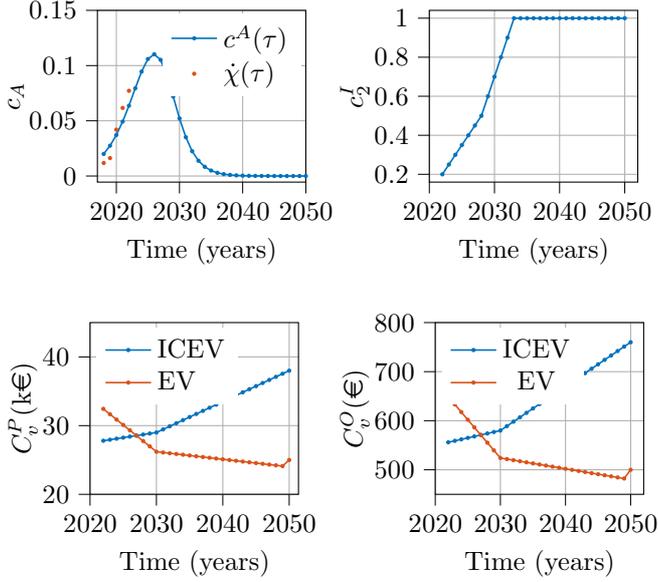
\begin{figure}[h!]
    \centering
    \begin{subfigure}[t]{0.48\columnwidth}
        \centering
        \pgfplotsset{ylabel shift = -.5em} 
\definecolor{mycolor1}{rgb}{0.00000,0.44700,0.74100}%
\definecolor{mycolor2}{rgb}{0.85000,0.32500,0.09800}%
\definecolor{mycolor3}{rgb}{0.92900,0.69400,0.12500}%
\definecolor{mycolor4}{rgb}{0.49400,0.18400,0.55600}%

\begin{tikzpicture}

\definecolor{darkgray176}{RGB}{176,176,176}
\definecolor{lightgray204}{RGB}{204,204,204}
\definecolor{steelblue31119180}{RGB}{31,119,180}

\begin{axis}[
/pgf/number format/1000 sep={},
scale = 0.4,
legend cell align={left},
legend style={fill opacity=0, draw opacity=1, text opacity=1, draw=none},
tick align=outside,
tick pos=left,
x grid style={darkgray176},
xlabel={Time (years)},
xmajorgrids,
xmin=2017, xmax=2050,
xtick style={color=black},
y grid style={darkgray176},
ylabel={ $c_A$},
ymajorgrids,
ymin=-0.0055125, ymax=0.15,
ytick style={color=black},
yticklabel style={
        /pgf/number format/fixed,
        /pgf/number format/precision=3
},
scaled y ticks=false
]
\addplot [semithick, mycolor1, mark=*, mark size=.5, mark options={solid}]
table {%
2018 0.02
2019 0.02744
2020 0.03713
2021 0.04927
2022 0.06369
2023 0.07946
2024 0.09457
2025 0.10597
2026 0.11025
2027 0.10516
2028 0.09125
2029 0.07201
2030 0.05212
2031 0.03514
2032 0.02247
2033 0.01384
2034 0.00833
2035 0.00494
2036 0.0029
2037 0.0017
2038 0.00099
2039 0.00057
2040 0.00033
2041 0.00019
2042 0.00011
2043 7e-05
2044 4e-05
2045 2e-05
2046 1e-05
2047 1e-05
2048 0
2049 0
2050 0
};
\addlegendentry{$c^A(\tau)$}
\addplot [semithick, mycolor2, mark=*, mark size=.5, mark options={solid}, only marks]
table {%
2018 0.0118069936061286
2019 0.0162375472072179
2020 0.041949571268378
2021 0.0616447665732201
2022 0.0771788740705876
};
\addlegendentry{$\dot{\chi}(\tau)$}
\end{axis}

\end{tikzpicture}
    \end{subfigure}
    \begin{subfigure}[t]{0.48\columnwidth}
        \centering
        \pgfplotsset{ylabel shift = -.5em} 
\definecolor{mycolor1}{rgb}{0.00000,0.44700,0.74100}%
\definecolor{mycolor2}{rgb}{0.85000,0.32500,0.09800}%
\definecolor{mycolor3}{rgb}{0.92900,0.69400,0.12500}%
\definecolor{mycolor4}{rgb}{0.49400,0.18400,0.55600}%

\begin{tikzpicture}

\definecolor{darkgray176}{RGB}{176,176,176}
\definecolor{steelblue31119180}{RGB}{31,119,180}

\begin{axis}[
scale=0.4,
/pgf/number format/1000 sep={},
tick align=outside,
tick pos=left,
x grid style={darkgray176},
xlabel={Time (years)},
xmajorgrids,
xmin=2020, xmax=2052,
xtick style={color=black},
y grid style={darkgray176},
ylabel={$c_2^I$},
ymajorgrids,
ymin=0.16, ymax=1.04,
ytick style={color=black}
]
\addplot [semithick, mycolor1, mark=*, mark size=.5, mark options={solid}]
table {%
2022 0.2
2023 0.25
2024 0.3
2025 0.35
2026 0.4
2027 0.45
2028 0.5
2029 0.6
2030 0.7
2031 0.8
2032 0.9
2033 1
2034 1
2035 1
2036 1
2037 1
2038 1
2039 1
2040 1
2041 1
2042 1
2043 1
2044 1
2045 1
2046 1
2047 1
2048 1
2049 1
2050 1
};
\end{axis}

\end{tikzpicture}
    \end{subfigure}
    
    \begin{subfigure}[t]{0.48\columnwidth}
        \centering
        \pgfplotsset{ylabel shift = -.5em} 
\definecolor{mycolor1}{rgb}{0.00000,0.44700,0.74100}%
\definecolor{mycolor2}{rgb}{0.85000,0.32500,0.09800}%
\definecolor{mycolor3}{rgb}{0.92900,0.69400,0.12500}%
\definecolor{mycolor4}{rgb}{0.49400,0.18400,0.55600}%

\begin{tikzpicture}

\definecolor{darkgray176}{RGB}{176,176,176}
\definecolor{darkorange25512714}{RGB}{255,127,14}
\definecolor{steelblue31119180}{RGB}{31,119,180}

\begin{axis}[
/pgf/number format/1000 sep={},
scale = 0.4,
tick align=outside,
tick pos=left,
x grid style={darkgray176},
xlabel={},
xmajorgrids,
xmin=2020, xmax=2051.4,
xtick style={color=black},
xlabel={Time (years)},
y grid style={darkgray176},
ylabel=
{$C_v^P$(k\euro)},
ymajorgrids,
ymin=20.000, ymax=45.000,
ytick style={color=black},
legend pos=north west, 
legend cell align={left},
legend style={fill opacity=0, draw opacity=1, text opacity=1, draw=none},
yticklabel style={
        /pgf/number format/fixed,
        /pgf/number format/precision=4
},
scaled y ticks=false
]
\addplot [semithick, mycolor1, mark=*, mark size=.5, mark options={solid}]
table {%
2022 27.800
2023 27.950
2024 28.100
2025 28.250
2026 28.400
2027 28.550
2028 28.700
2029 28.850
2030 29.000
2031 29.450
2032 29.900
2033 30.350
2034 30.800
2035 31.250
2036 31.700
2037 32.150
2038 32.600
2039 33.050
2040 33.500
2041 33.950
2042 34.400
2043 34.850
2044 35.300
2045 35.750
2046 36.200
2047 36.650
2048 37.100
2049 37.550
2050 38.000
};
\addlegendentry{ICEV}
\addplot [semithick, mycolor2, mark=*, mark size=.5, mark options={solid}]
table {%
2022 32.440
2023 31.660
2024 30.880
2025 30.100
2026 29.320
2027 28.540
2028 27.760
2029 26.980
2030 26.200
2031 26.090
2032 25.980
2033 25.870
2034 25.760
2035 25.650
2036 25.540
2037 25.430
2038 25.320
2039 25.210
2040 25.100
2041 24.990
2042 24.880
2043 24.770
2044 24.660
2045 24.550
2046 24.440
2047 24.330
2048 24.220
2049 24.110
2050 25.000
};
\addlegendentry{EV}
\end{axis}

\end{tikzpicture}
    \end{subfigure}
    \begin{subfigure}[t]{0.48\columnwidth}
        \centering
        \pgfplotsset{ylabel shift = -.5em} 
\definecolor{mycolor1}{rgb}{0.00000,0.44700,0.74100}%
\definecolor{mycolor2}{rgb}{0.85000,0.32500,0.09800}%
\definecolor{mycolor3}{rgb}{0.92900,0.69400,0.12500}%
\definecolor{mycolor4}{rgb}{0.49400,0.18400,0.55600}%

\begin{tikzpicture}

\definecolor{darkgray176}{RGB}{176,176,176}
\definecolor{darkorange25512714}{RGB}{255,127,14}
\definecolor{steelblue31119180}{RGB}{31,119,180}


\begin{axis}[
/pgf/number format/1000 sep={},
scale = 0.4,
tick pos=left,
tick align=outside,
x grid style={darkgray176},
xmin=2020, xmax=2051.4,
xtick pos=left,
xlabel={Time (years)},
xmajorgrids,
xmin=2020, xmax=2052,
xtick style={color=black},
y grid style={darkgray176},
ymajorgrids,
ylabel=
{$C_v^O$(\euro)},
ymin=450, ymax=800,
ytick style={color=black},
legend pos=north west, 
legend style={fill opacity=0, draw opacity=1, text opacity=1, draw=none},
]
\addplot [semithick, mycolor1, mark=*, mark size=.5, mark options={solid}]
table {%
2022 556
2023 559
2024 562
2025 565
2026 568
2027 571
2028 574
2029 577
2030 580
2031 589
2032 598
2033 607
2034 616
2035 625
2036 634
2037 643
2038 652
2039 661
2040 670
2041 679
2042 688
2043 697
2044 706
2045 715
2046 724
2047 733
2048 742
2049 751
2050 760
};
\addlegendentry{ICEV}
\addplot [semithick, mycolor2, mark=*, mark size=.5, mark options={solid}]
table {%
2022 648.8
2023 633.2
2024 617.6
2025 602
2026 586.4
2027 570.8
2028 555.2
2029 539.6
2030 524
2031 521.8
2032 519.6
2033 517.4
2034 515.2
2035 513
2036 510.8
2037 508.6
2038 506.4
2039 504.2
2040 502
2041 499.8
2042 497.6
2043 495.4
2044 493.2
2045 491
2046 488.8
2047 486.6
2048 484.4
2049 482.2
2050 500
};
\addlegendentry{EV}
\end{axis}

\end{tikzpicture}
    \end{subfigure}
%
\caption{Determinants used in Logit model} \label{fig:Exo}
\end{figure}

The adoption coefficient $c^A(t)$ is the solution of the Bass model (normalised to market-share) 
\begin{equation}\label{eqn:Bass}
\begin{aligned}
    c^A(\tau) = \frac{d}{d\tau}\chi(\tau) = (p + q \chi(\tau))(1-\chi(\tau)) 
\end{aligned}
\end{equation}
where 
$p$ and $q$ are the coefficients of innovation and imitation, respectively. The values of $p$ and $q$ are tuned to align with the yearly EV sales from 2018 to 2022, as provided by data source 4. Figure~\ref{fig:Exo} and Table~\ref{tab:logit} shows the different determinants and parameters, used in  (\ref{eqn:U}), respectively.

\begin{table}[h!]
    \centering
    \caption{List of model parameters.}   
    \label{tab:logit}
    \begin{tabular}{ccc}
        Attributes & ICEV & EV \\ \hline
        Purchase Cost ($p^P$) & -0.3 & -0.3\\
        Operating Cost ($p^O$) & -0.15 & -0.15\\
        Infrastructure Cost ($p^I$) & - & -0.3\\
        $\mu$ & \multicolumn{2}{c}{6.75} \\ 
        $p$ & \multicolumn{2}{c}{0.02} \\ 
        $q$ & \multicolumn{2}{c}{0.4} \\ \hline
    \end{tabular}
    \label{tab:param}
\end{table}

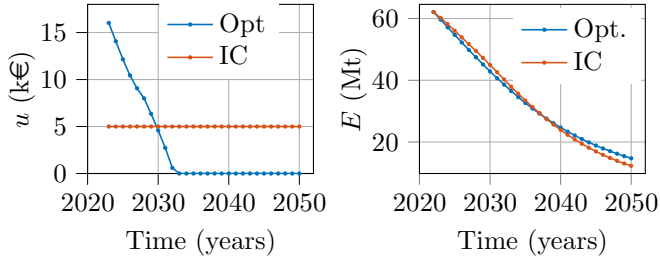
\begin{figure}[h!]
    \centering
    \begin{subfigure}[t]{.48\columnwidth}
        \centering
         \pgfplotsset{ylabel shift = -.5em} 
\definecolor{mycolor1}{rgb}{0.00000,0.44700,0.74100}%
\definecolor{mycolor2}{rgb}{0.85000,0.32500,0.09800}%
\definecolor{mycolor3}{rgb}{0.92900,0.69400,0.12500}%
\definecolor{mycolor4}{rgb}{0.49400,0.18400,0.55600}%

\begin{tikzpicture}

\definecolor{darkgray176}{RGB}{176,176,176}
\definecolor{lightgray204}{RGB}{204,204,204}
\definecolor{steelblue31119180}{RGB}{31,119,180}

\begin{axis}[
/pgf/number format/1000 sep={},
width=4.521in,
height=3.566in,
at={(0.758in,0.481in)},
scale=0.3,
legend cell align={left},
legend style={
  fill opacity=0,
  draw opacity=1,
  text opacity=1,
  at= {(axis cs: 2030,10)},
  anchor=south west,
  draw=none
},
tick align=outside,
tick pos=left,
x grid style={darkgray176},
xlabel={Time (years)},
xmajorgrids,
xmin=2020, xmax=2052,
xtick style={color=black},
y grid style={darkgray176},
ylabel={$u$ (k\euro) },
ymajorgrids,
ymin=0, ymax=18,
ytick style={color=black}
]
\addplot [semithick, mycolor1, mark=*, mark size=.5, mark options={solid}]
table {%

2023 16.027071352215
2024 14.06938781574
2025 12.15172890585
2026 10.44599065068
2027 9.072187285755
2028 8.0056819002
2029 6.354953454225
2030 4.5890303616
2031 2.732734600561
2032 0.604945902374
2033 5.89397748981e-09
2034 1.143075510108e-08
2035 2.031528948005e-10
2036 6.3794688345e-09
2037 5.25227425541e-09
2038 8.80376114192e-11
2039 3.65187944961e-09
2040 1.70971813564e-09
2041 2.858319209011e-09
2042 2.092646033556e-09
2043 5.84913549451e-10
2044 2.146455329164e-09
2045 2.08022528754e-09
2046 1.692750388456e-10
2047 1.276257252809e-09
2048 7.0476540015e-10
2049 2.293180393302e-10
2050 9.859974642e-11

};
\addlegendentry{Opt}

\addplot [semithick, mycolor2, mark=*, mark size=.5, mark options={solid}]
table {%
2023 5
2024 5
2025 5
2026 5
2027 5
2028 5
2029 5
2030 5
2031 5
2032 5
2033 5
2034 5
2035 5
2036 5
2037 5
2038 5
2039 5
2040 5
2041 5
2042 5
2043 5
2044 5
2045 5
2046 5
2047 5
2048 5
2049 5
2050 5
};
\addlegendentry{IC}

\end{axis}

\end{tikzpicture}
    \end{subfigure}
    \begin{subfigure}[t]{.48\columnwidth}
        \centering
         \pgfplotsset{ylabel shift = -.5em} 
\definecolor{mycolor1}{rgb}{0.00000,0.44700,0.74100}%
\definecolor{mycolor2}{rgb}{0.85000,0.32500,0.09800}%
\definecolor{mycolor3}{rgb}{0.92900,0.69400,0.12500}%
\definecolor{mycolor4}{rgb}{0.49400,0.18400,0.55600}%

\begin{tikzpicture}

\definecolor{darkgray176}{RGB}{176,176,176}
\definecolor{lightgray204}{RGB}{204,204,204}
\definecolor{steelblue31119180}{RGB}{31,119,180}

\begin{axis}[
/pgf/number format/1000 sep={},
width=4.521in,
height=3.566in,
at={(0.758in,0.481in)},
scale=0.3,
legend cell align={left},
legend style={fill opacity=0, draw opacity=1, text opacity=1, draw=none},
tick align=outside,
tick pos=left,
x grid style={darkgray176},
xlabel={Time (years)},
xmajorgrids,
xmin=2020, xmax=2052,
y grid style={darkgray176},
ylabel={$E$ (Mt)},
ymajorgrids,
ymin=9.79146907119933, ymax=64.5880943730381,
]
\addplot [semithick, mycolor1, mark=*, mark options={dashed, mycolor1}, mark size=.5, mark options={solid}]
table {%
2022 62.0973386775
2023 59.6369273092608
2024 57.1515060019351
2025 54.6786129429341
2026 52.2284823100033
2027 49.8241624410995
2028 47.4653029682664
2029 45.1409349468949
2030 42.8674930762398
2031 40.6816046421603
2032 38.5899245895381
2033 36.5040089862257
2034 34.5256276910336
2035 32.6470704667808
2036 30.865066296883
2037 29.1971991201004
2038 27.5899303856052
2039 26.0877697603396
2040 24.6833631878252
2041 23.3707392667316
2042 22.1288935954514
2043 20.9668524921619
2044 19.8804295251482
2045 18.8665690866638
2046 17.9177302468571
2047 17.0357915775675
2048 16.2045649240206
2049 15.4275600754105
2050 14.7298093987136
};
\addlegendentry{Opt.}

\addplot [semithick,color=mycolor2, mark=*, mark options={dashed, mycolor2}, mark size=.5]
  table[row sep=crcr]{%
2022	62.0973386775\\
2023	60.1836066576528\\
2024	58.1593466433033\\
2025	56.0574852427096\\
2026	53.8927857845427\\
2027	51.7029605179969\\
2028	49.5027172631407\\
2029	47.2461952797078\\
2030	44.9385068147861\\
2031	42.6121506008471\\
2032	40.2665216203798\\
2033	37.8983092627524\\
2034	35.6270367242917\\
2035	33.446161445072\\
2036	31.3547818773782\\
2037	29.3718667035712\\
2038	27.4509062460178\\
2039	25.6416200529808\\
2040	23.9416930704035\\
2041	22.347770786495\\
2042	20.8449561356326\\
2043	19.4437224049606\\
2044	18.1443956814805\\
2045	16.9438462199866\\
2046	15.8372562880573\\
2047	14.8249199677192\\
2048	13.8917475730213\\
2049	13.0384358223612\\
2050	12.2822247684378\\
};
\addlegendentry{IC}


\end{axis}

\end{tikzpicture}
    \end{subfigure}
    \caption{Incentive (left) and Emission (right) Profile} \label{fig:opt}
\end{figure}

\begin{figure}[h!]
    \centering
    \begin{subfigure}[t]{\columnwidth}
        \centering
\begin{tikzpicture}

\definecolor{mycolor1}{rgb}{0.00000,0.44700,0.74100}%
\definecolor{mycolor2}{rgb}{0.85000,0.32500,0.09800}%
\definecolor{mycolor3}{rgb}{0.92900,0.69400,0.12500}%
\definecolor{mycolor4}{rgb}{0.49400,0.18400,0.55600}%

\begin{axis}[
/pgf/number format/1000 sep={},
width=4.521in,
height=3.566in,
at={(0.758in,0.481in)},
scale=0.4,
legend cell align={left},
legend style={
  fill opacity=0,
  draw opacity=1,
  text opacity=1,
  at={(1.1,0.97)},
  anchor=north west,
  draw=none
},
tick align=outside,
tick pos=left,
x grid style={black},
xlabel={Time (years)},
xmajorgrids,
xmin=2020, xmax=2052,
xtick style={color=black},
y grid style={black},
ylabel={$N_v$ (Mvehicles)},
ymajorgrids,
ymin=0, ymax=2,
ytick style={color=black}
]

\addplot [semithick, mycolor1, mark=*, mark size=.5, mark options={solid}]
table {%

2023 1.01517927904732
2024 1.06467183121248
2025 1.11519865799253
2026 1.1557135713313
2027 1.18471142318081
2028 1.20203520853632
2029 1.24291730044253
2030 1.2831198606974
2031 1.30174932808686
2032 1.32345233175703
2033 1.4154064917187
2034 1.43233137434477
2035 1.45625342319946
2036 1.48512134800008
2037 1.51724599382448
2038 1.55123387348811
2039 1.58592484972935
2040 1.62039571128848
2041 1.65384700243842
2042 1.68571720856938
2043 1.71560340930812
2044 1.74330315840776
2045 1.76868313381304
2046 1.79172890567911
2047 1.81252999716017
2048 1.83127810029943
2049 1.84824414494137
2050 1.81754743365601
};
\addlegendentry{EV (Opt.)}

\addplot [semithick, mycolor3, mark=*, mark size=.5, mark options={solid}]
table {%

2023 1.70315855846634
2024 1.57921548136421
2025 1.48032603311622
2026 1.39129158658586
2027 1.30967557788741
2028 1.23190130990965
2029 1.13478753837855
2030 1.04868262843489
2031 0.994418447912606
2032 0.947166516038613
2033 0.839136296413857
2034 0.814580815062154
2035 0.790137264871424
2036 0.766330126687121
2037 0.743304927651748
2038 0.720980456038655
2039 0.699177932177081
2040 0.677722776854636
2041 0.656449011236205
2042 0.635271536293914
2043 0.614156973973312
2044 0.593136030210035
2045 0.572251211560599
2046 0.551571730143368
2047 0.531183647157362
2048 0.511185095232278
2049 0.491676353888846
2050 0.518879603694141

};
\addlegendentry{ICEV (Opt.)}

\addplot [semithick, mycolor2, mark=*, mark size=.5, mark options={solid}]
table {%
2023 0.62580652378235
2024 0.733020415911391
2025 0.845749378103449
2026 0.946376442791526
2027 1.02582653254529
2028 1.08335491076081
2029 1.18873546521644
2030 1.29958796087198
2031 1.3905438807877
2032 1.48972475109184
2033 1.59262823621004
2034 1.6063756014762
2035 1.62690518534374
2036 1.65228975273203
2037 1.68088941469848
2038 1.71131195105727
2039 1.74237556984323
2040 1.77313400618992
2041 1.80276439570704
2042 1.83069824918656
2043 1.85653620962151
2044 1.88009327621429
2045 1.90125579787058
2046 1.92003479917163
2047 1.93654886404368
2048 1.9510212746768
2049 1.96375474999276
2050 1.9357186557451
};
\addlegendentry{EV (IC)}

\end{axis}

\end{tikzpicture}
    \end{subfigure}
    \hspace{1cm}
    \begin{subfigure}[t]{\columnwidth}
        \centering
\definecolor{mycolor1}{rgb}{0.00000,0.44700,0.74100}%
\definecolor{mycolor2}{rgb}{0.85000,0.32500,0.09800}%
\definecolor{mycolor3}{rgb}{0.92900,0.69400,0.12500}%
\definecolor{mycolor4}{rgb}{0.49400,0.18400,0.55600}%

\begin{tikzpicture}

\definecolor{darkgray176}{RGB}{176,176,176}
\definecolor{lightgray204}{RGB}{204,204,204}
\definecolor{steelblue31119180}{RGB}{31,119,180}

\begin{axis}[
/pgf/number format/1000 sep={},
width=4.521in,
height=3.566in,
at={(0.758in,0.481in)},
scale=0.4,
legend cell align={left},
legend style={fill opacity=0.2, draw opacity=1, text opacity=1, draw=none,
at={(1.1,0.97)},
 anchor=north west},
tick align=outside,
tick pos=left,
x grid style={darkgray176},
xlabel={Time (years)},
xmajorgrids,
xmin=2020, xmax=2052,
xtick style={color=black},
y grid style={darkgray176},
ylabel={$S_v$ (Mvehicles)},
ymajorgrids,
ymin=0, ymax= 40,
ytick style={color=black}
]
\addplot [semithick, mycolor1, mark=*, mark size=.5, mark options={solid}]
table {%

2022 0.276883
2023 1.29135486794373
2024 2.35490365133818
2025 3.46841519676797
2026 4.62167958963476
2027 5.80261408076466
2028 6.99862304917957
2029 8.22592828929042
2030 9.47330873035929
2031 10.7086343782447
2032 11.9247898406515
2033 13.182615207678
2034 14.3988247635827
2035 15.5733888339143
2036 16.7055196332013
2037 17.7942620455312
2038 18.8388496202931
2039 19.8383784348856
2040 20.7922292117582
2041 21.7002550093439
2042 22.562786288021
2043 23.3805539520607
2044 24.154661719457
2045 24.8864443852979
2046 25.5774007359189
2047 26.2291441522406
2048 26.8433829691061
2049 27.4219125398571
2050 27.9204094245999
};
\addlegendentry{EV (Opt.)}

\addplot [semithick, mycolor3, mark=*, mark size=.5, mark options={solid}]
table {%

2022 35.519509
2023 34.7825329307659
2024 33.9617324389014
2025 33.083521296982
2026 32.158109419604
2027 31.1975795559451
2028 30.2145273269816
2029 29.1927309383006
2030 28.1434114606382
2031 27.0986988881433
2032 26.0657086131055
2033 24.9836005454302
2034 23.9356604008552
2035 22.9219178538336
2036 21.943160689837
2037 21.0003440247777
2038 20.0942343092667
2039 19.2257354659038
2040 18.3954667722428
2041 17.6035751698473
2042 16.8497301983353
2043 16.1332009534464
2044 15.4528837171812
2045 14.80744369445
2046 14.195382098919
2047 13.6150855496677
2048 13.0648457118528
2049 12.5428672321329
2050 12.0934735503995
};
\addlegendentry{ICEV (Opt.)}
\addplot [color=mycolor2, mark=*, mark options={solid, mycolor2}, mark size=.5]
  table[row sep=crcr]{%
2022	0.276883\\
2023	0.90198211267876\\
2024	1.63387948077211\\
2025	2.47794174631282\\
2026	3.42186901063984\\
2027	4.44391861113421\\
2028	5.52124728177362\\
2029	6.69710325793011\\
2030	7.97000666968697\\
2031	9.31230240582858\\
2032	10.7239404782539\\
2033	12.2002755811868\\
2034	13.6441437964197\\
2035	15.0546221112408\\
2036	16.4291370157325\\
2037	17.7640184007189\\
2038	19.0548684113387\\
2039	20.2967919404005\\
2040	21.4852832249128\\
2041	22.6165959604329\\
2042	23.6878800947605\\
2043	24.6972043125328\\
2044	25.6436017441817\\
2045	26.5269593070818\\
2046	27.3479553979408\\
2047	28.1079858433601\\
2048	28.8090953418799\\
2049	29.4539002713503\\
2050	30.0061194771407\\
};
\addlegendentry{EV (IC)}

\end{axis}

\end{tikzpicture}
    \end{subfigure}
    \caption{Vehicle Sales (top) and Vehicle Stock (bottom) }\label{fig:salesStock}
\end{figure}
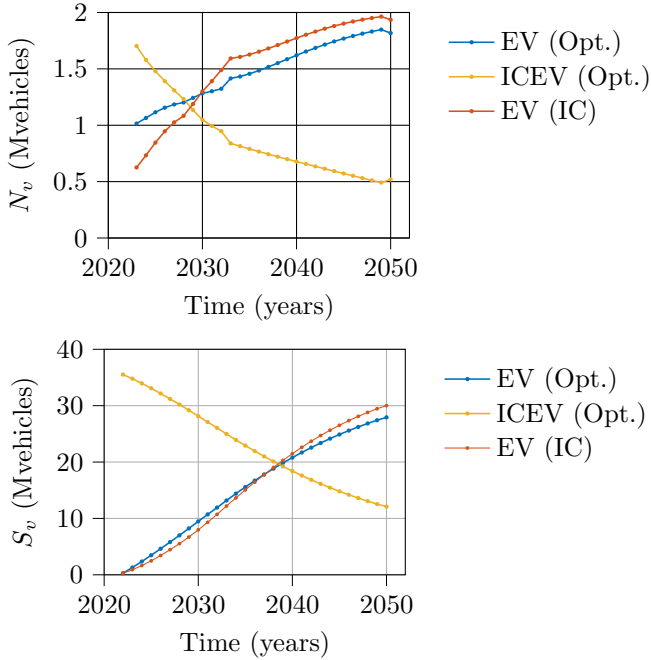

\subsection{Results}

The emissions, and vehicle sales and stock curves, for IC and optimal scenario, are shown in Fig.~\ref{fig:opt}~(right), and Fig.~\ref{fig:salesStock}, respectively. Clearly, in both scenarios, the ICEV stock ($v=1$) decreases while the EV stock ($v=2$) increases with time, both exhibiting an S-shaped curve suggesting variable rate. The curve of CO$_2$ emissions, Fig.~\ref{fig:opt}~(right), is proportional to that of $S_1$ and decreases by more than three times with respect to 2022. The IC scenario forecasts 12.3~Mt of CO$_2$ in 2050. Correspondingly, $\overline{E}$ is set to this value. 

The incentive law shown in  Fig.~\ref{fig:opt}~(left) exhibits a very variable behavior, being null until a certain year, to rise up by the end of the period. Intuitively, such a behavior is optimal, within the assumptions of the model, in that it incentivises late adopters and encourages EV purchase during a period when ICE performance is improving. 
and further reduce the emissions while minimizing the total budget. This effect is visible in the curves of EV sales and stock (Fig.~\ref{fig:salesStock}) and emissions (Fig.~\ref{fig:opt}~(right)). Compared to the IC law,  EV sales are lower in the initial years and begin to increase around midpoint. Consequently, the optimal emission curve decreases gradually in the early years and exhibits a sharper decrease towards the end. The final CO$_2$ emission for both scenarios are equal as the  constraint is defined only at $T$. However, the cumulative CO$_2$ emissions of the optimal scenario is higher than the IC scenario. To address this, an integral constraint on emissions can be defined in future work.  
Regarding the EV stock, it increases gradually early-on but shows a sharper increase towards the end as result of the optimal incentive. 

Overall, we obtain a total expenditure $I(T)= 196.2$ G\euro, that is, around 30\% reduction with respect to the IC reference scenario, see Table~\ref{tab:res}.    

\subsection{Alternative scenarios}\label{sec:alter}

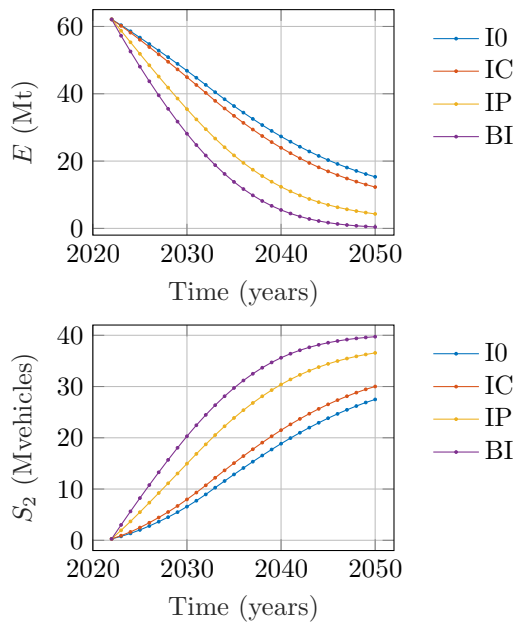
\begin{figure}[h!] 
    \centering
    \begin{subfigure}[t]{\columnwidth}
        \centering
%
%
\definecolor{mycolor1}{rgb}{0.00000,0.44700,0.74100}%
\definecolor{mycolor2}{rgb}{0.85000,0.32500,0.09800}%
\definecolor{mycolor3}{rgb}{0.92900,0.69400,0.12500}%
\definecolor{mycolor4}{rgb}{0.49400,0.18400,0.55600}%
\begin{tikzpicture}

\begin{axis}[%
/pgf/number format/1000 sep={},
width=4.521in,
height=3.566in,
at={(0.758in,0.481in)},
scale =.4,
xmin=2020,
xmax=2052,
xlabel style={font=\color{white!15!black}},
xlabel={Time (years)},
ymin=-2,
ymax=65,
ylabel style={font=\color{white!15!black}},
ylabel={$E$ (Mt)},
axis background/.style={fill=white},
xmajorgrids,
ymajorgrids,
legend style={at={(1.1,.97)}, anchor=north west, legend cell align=left, align=left, draw=none, fill opacity=0.8}
]
\addplot [color=mycolor1, mark=*, mark options={solid, mycolor1}, mark size=.5]
  table[row sep=crcr]{%
2022	62.0973386775\\
2023	60.3772061729652\\
2024	58.5655250320406\\
2025	56.6930285796768\\
2026	54.7698219810999\\
2027	52.8296600717508\\
2028	50.8844656595381\\
2029	48.8888122329172\\
2030	46.840855906383\\
2031	44.760913732486\\
2032	42.6405497339573\\
2033	40.468819463666\\
2034	38.375787677183\\
2035	36.3531099403197\\
2036	34.3985618075106\\
2037	32.532486521514\\
2038	30.7057871997411\\
2039	28.9683729218248\\
2040	27.3185351239776\\
2041	25.7558635321482\\
2042	24.2646197819138\\
2043	22.8586192926125\\
2044	21.5378434203587\\
2045	20.3026815623895\\
2046	19.1482585401559\\
2047	18.0783267996862\\
2048	17.0778073973195\\
2049	16.1506154509676\\
2050	15.3215713193231\\
};
\addlegendentry{I0}

\addplot [color=mycolor2, mark=*, mark options={solid, mycolor2}, mark size=.5]
  table[row sep=crcr]{%
2022	62.0973386775\\
2023	60.1836066576528\\
2024	58.1593466433033\\
2025	56.0574852427096\\
2026	53.8927857845427\\
2027	51.7029605179969\\
2028	49.5027172631407\\
2029	47.2461952797078\\
2030	44.9385068147861\\
2031	42.6121506008471\\
2032	40.2665216203798\\
2033	37.8983092627524\\
2034	35.6270367242917\\
2035	33.446161445072\\
2036	31.3547818773782\\
2037	29.3718667035712\\
2038	27.4509062460178\\
2039	25.6416200529808\\
2040	23.9416930704035\\
2041	22.347770786495\\
2042	20.8449561356326\\
2043	19.4437224049606\\
2044	18.1443956814805\\
2045	16.9438462199866\\
2046	15.8372562880573\\
2047	14.8249199677192\\
2048	13.8917475730213\\
2049	13.0384358223612\\
2050	12.2822247684378\\
};
\addlegendentry{IC}

\addplot [color=mycolor3, mark=*, mark options={solid, mycolor3}, mark size=.5]
  table[row sep=crcr]{%
2022	62.0973386775\\
2023	58.7223564746349\\
2024	55.3012184625432\\
2025	51.8747599714611\\
2026	48.4665169250132\\
2027	45.1128889718546\\
2028	41.824636303344\\
2029	38.588923706835\\
2030	35.4343880572894\\
2031	32.3903154323214\\
2032	29.4742336848917\\
2033	26.6935893006537\\
2034	24.1001622652403\\
2035	21.6911997086484\\
2036	19.466136748236\\
2037	17.4358212160008\\
2038	15.5542730620576\\
2039	13.8638097030724\\
2040	12.3533559832307\\
2041	11.0053001882722\\
2042	9.79984113216183\\
2043	8.73318347113825\\
2044	7.80070907455796\\
2045	6.98622610477781\\
2046	6.28097052761465\\
2047	5.67424371690844\\
2048	5.14866273304823\\
2049	4.69634039844939\\
2050	4.30014740913832\\
};
\addlegendentry{IP}

\addplot [color=mycolor4, mark=*, mark options={solid, mycolor4}, mark size=.5]
  table[row sep=crcr]{%
2022	62.0973386775\\
2023	57.245692693174\\
2024	52.5643722590114\\
2025	48.0530702524094\\
2026	43.7059135411688\\
2027	39.5335316421169\\
2028	35.5282360685896\\
2029	31.7193255315308\\
2030	28.1293588279806\\
2031	24.7717666282476\\
2032	21.6599659491601\\
2033	18.7950313188431\\
2034	16.1808434546874\\
2035	13.8145286760166\\
2036	11.6928169436272\\
2037	9.8194005597368\\
2038	8.14413720479249\\
2039	6.70151432013578\\
2040	5.47205851709268\\
2041	4.4276239453184\\
2042	3.54232659076817\\
2043	2.8027270151109\\
2044	2.19929010791141\\
2045	1.7079486749557\\
2046	1.31683557132408\\
2047	1.00957371169903\\
2048	0.767686614436828\\
2049	0.579754571869599\\
2050	0.419126668371727\\
};
\addlegendentry{BI}

\end{axis}

\end{tikzpicture}%
    \end{subfigure}
    \hspace{1cm}
    \begin{subfigure}[t]{\columnwidth}
        \centering
%
%
\definecolor{mycolor1}{rgb}{0.00000,0.44700,0.74100}%
\definecolor{mycolor2}{rgb}{0.85000,0.32500,0.09800}%
\definecolor{mycolor3}{rgb}{0.92900,0.69400,0.12500}%
\definecolor{mycolor4}{rgb}{0.49400,0.18400,0.55600}%
\begin{tikzpicture}

\begin{axis}[%
/pgf/number format/1000 sep={},
width=4.521in,
height=3.566in,
at={(0.758in,0.481in)},
scale =.4,
xmin=2020,
xmax=2052,
xlabel style={font=\color{white!15!black}},
xlabel={Time (years)},
ymin=-2,
ymax=42,
ylabel style={font=\color{white!15!black}},
ylabel={$S_2$ (Mvehicles)},
axis background/.style={fill=white},
xmajorgrids,
ymajorgrids,
legend style={at={(1.1,0.97)}, 
anchor=north west, 
legend cell align=left, 
align=left, 
draw=none, 
fill opacity=0.8}
]
\addplot [color=mycolor1, mark=*, mark options={solid, mycolor1},, mark size=.5]
  table[row sep=crcr]{%
2022	0.276883\\
2023	0.764090720006111\\
2024	1.34310863939405\\
2025	2.02060231508006\\
2026	2.78741710446031\\
2027	3.62453088480891\\
2028	4.51102615279705\\
2029	5.48964332250078\\
2030	6.56396393107977\\
2031	7.71554823951511\\
2032	8.95042866333868\\
2033	10.2700569397255\\
2034	11.5690713659481\\
2035	12.848016378594\\
2036	14.105449518351\\
2037	15.3384870919928\\
2038	16.5431582723137\\
2039	17.7146193606683\\
2040	18.8481157382805\\
2041	19.9394168507441\\
2042	20.984976630418\\
2043	21.9819912260398\\
2044	22.9284789709706\\
2045	23.8232150217884\\
2046	24.6657077362019\\
2047	25.4561605230879\\
2048	26.1954340401284\\
2049	26.8849951828921\\
2050	27.4806489361322\\
};
\addlegendentry{I0}

\addplot [color=mycolor2, mark=*, mark options={solid, mycolor2}, mark size=.5]
  table[row sep=crcr]{%
2022	0.276883\\
2023	0.90198211267876\\
2024	1.63387948077211\\
2025	2.47794174631282\\
2026	3.42186901063984\\
2027	4.44391861113421\\
2028	5.52124728177362\\
2029	6.69710325793011\\
2030	7.97000666968697\\
2031	9.31230240582858\\
2032	10.7239404782539\\
2033	12.2002755811868\\
2034	13.6441437964197\\
2035	15.0546221112408\\
2036	16.4291370157325\\
2037	17.7640184007189\\
2038	19.0548684113387\\
2039	20.2967919404005\\
2040	21.4852832249128\\
2041	22.6165959604329\\
2042	23.6878800947605\\
2043	24.6972043125328\\
2044	25.6436017441817\\
2045	26.5269593070818\\
2046	27.3479553979408\\
2047	28.1079858433601\\
2048	28.8090953418799\\
2049	29.4539002713503\\
2050	30.0061194771407\\
};
\addlegendentry{IC}

\addplot [color=mycolor3, mark=*, mark options={solid, mycolor3}, mark size=.5]
  table[row sep=crcr]{%
2022	0.276883\\
2023	1.94275859630977\\
2024	3.67924280204581\\
2025	5.48524921463686\\
2026	7.34119964602009\\
2027	9.22532531067205\\
2028	11.1167276166683\\
2029	13.033159738462\\
2030	14.9544919695955\\
2031	16.8546197236402\\
2032	18.7185467046846\\
2033	20.532052113359\\
2034	22.2505828732729\\
2035	23.8698426699347\\
2036	25.3864047645767\\
2037	26.7979073407248\\
2038	28.1031481145103\\
2039	29.3020929008776\\
2040	30.396113327194\\
2041	31.3878861499674\\
2042	32.2812878614017\\
2043	33.0811692224301\\
2044	33.7931850646948\\
2045	34.4234771791942\\
2046	34.9784588380751\\
2047	35.4646313774455\\
2048	35.8884526528261\\
2049	36.2562371792307\\
2050	36.5649091272265\\
};
\addlegendentry{IP}

\addplot [color=mycolor4, mark=*, mark options={solid, mycolor4}, mark size=.5]
  table[row sep=crcr]{%
2022	0.276883\\
2023	2.99451342641007\\
2024	5.63727769116873\\
2025	8.23111526971474\\
2026	10.7756712491674\\
2027	13.2662813181847\\
2028	15.6941915965093\\
2029	18.0443318126144\\
2030	20.2998858464384\\
2031	22.4457899245698\\
2032	24.469308255942\\
2033	26.3603083848646\\
2034	28.1114522557907\\
2035	29.7182529965896\\
2036	31.1789829659319\\
2037	32.4945500232209\\
2038	33.6682873682154\\
2039	34.7056635346217\\
2040	35.613985795169\\
2041	36.4019224158842\\
2042	37.0791797602718\\
2043	37.6560781425844\\
2044	38.1432227102626\\
2045	38.5510590988991\\
2046	38.8895814344837\\
2047	39.1681168349862\\
2048	39.3952037070307\\
2049	39.5785361600049\\
2050	39.7248865118169\\
};
\addlegendentry{BI}

\end{axis}

\end{tikzpicture}%
    \end{subfigure}
    \caption{Reference Scenarios with the full model: CO$_2$ emissions (top) and EV stock (bottom) as a function of time.}\label{fig:refs}
\end{figure}

In addition to the optimal and IC scenarios, we analyze three different policies, namely,
\begin{itemize}
    \item No incentive (I0), $u(t) \equiv 0$ 
    \item Incentive covering the whole EV price (IP), $u(t) = C_2^P(t)$
    \item Ban of ICEV sales from $t_0$ (BI), $N_{2}(t) = N(t)$
\end{itemize}
The corresponding curves of CO$_2$ emissions and EV stock are shown in Fig.~\ref{fig:refs}. Clearly, the more stringent the policy, the faster increase of EV stock is observed, together with a decrease in yearly emissions. With the BI, the ICEV stock virtually empties by 2050 and consequently the CO$_2$ emissions vanish by that target year. The $E(T)$ and $I(T)$ values for the different policy scenarios are given in Table~\ref{tab:res}.

\begin{table}[h!]
    \centering
    \caption{Reference and optimal scenario}    \label{tab:res}
    \begin{tabular}{c|ccccc}
    \hline
              Output  & I0  & IC  & IP  & BI  & Optimal  \\
    \hline
       $E(T)$ (Mt) &   15.3 & 12.3 & 4.3 & 0.4 & 12.3 \\
         $I(T)$ (G\euro) & 0 & 215 & 1497 & - & 196.2 \\
         \hline
    \end{tabular}
\end{table}

\section{Conclusions and Future Work}

This study represents a first attempt to build a backcasting methodology to identify the optimal policy roadmaps in transport systems. 
The analysis focussed on a passenger car fleet subsystem, describing its evolution and associated emissions, with the monetary incentive to an EV purchase as the control input. The optimal incentive trajectory was derived by formulating an optimal control problem with the objective to minimize the state's budget while reaching a desired CO$_2$ target. A quantitative case study applied to Metropolitan France was performed to illustrate the backcasting approach. 

Further research can improve the backcasting paradigm in several ways. Refinements to the fleet model could include regional disaggregation within Metropolitan France, additional vehicle types (e.g., gasoline, diesel, hybrid) and transport modes (e.g., bikes, rail), and modeling second-hand vehicle exchanges. Mileage assumptions could be refined by accounting for variation by vehicle type and user profile, using zone-specific data (e.g., urban, rural) as in (\cite{ITF}). The survival rate, currently based only on natural obsolescence, could also incorporate factors like Low Emission Zones (LEZ), which can accelerate vehicle turnover. Additionally, the zonal choice in a LEZ implementation could be optimized.  Finally, to capture the changes in operating costs and vehicle ownership, the demand for passenger cars, treated here as exogenous, could be replaced with a demand model predicting vkm by mode and zone.

\section*{Acknowledgments}
This research benefited from state aid managed by the \textit{Agence Nationale de la Recherche (ANR)}, under France 2030, within the project FORBAC bearing the reference ANR-23-PEMO-0002. The authors would like to acknowledge the useful discussions with Dr. Benoit Cheze (IFPEN) and the partners of the FORBAC project. 

\bibliography{parc_aac}

\end{document}